\documentclass{amsart}
\usepackage{amsmath}
\usepackage{amsfonts}
\usepackage{amssymb}
\usepackage{graphicx}%

\usepackage{times}
\usepackage{bbold}
\usepackage{stmaryrd}
\usepackage{array}
\usepackage{dsfont}
\usepackage[T1]{fontenc}

\newtheorem{theorem}{Theorem}

\newtheorem*{GGRconjecture}{GGR Conjecture}
\newtheorem*{GGR+theorem}{GGR+ Theorem}

\newtheorem{lemma}[theorem]{Lemma}

\begin{document}
\title[A new proof of the GGR conjecture]{A new proof of the GGR conjecture}

\author{J. M. Ash}
\address{Department of Mathematics, DePaul University\\Chicago, IL 60614}
\email{mash@depaul.edu}
%\urladdr{http://www.depaul.edu/\symbol{126}mash/}
\author{S. Catoiu}
\address{Department of Mathematics, DePaul University\\Chicago, IL 60614}
\email{scatoiu@depaul.edu}
%\urladdr{http://www.depaul.edu/\symbol{126}scatoiu/}
\author{H. Fejzi\'{c}}
\address{Department of Mathematics, California State University, San Bernardino, CA 92407}
\email{hfejzic@csusb.edu}
\thanks{This paper is in final form and no
version of it will be submitted for publication elsewhere.}
\date{August 24, 2022}
\subjclass[2010]{Primary 26A24; Secondary 13F20; 15A03; 26A27.}
\keywords{GGR Conjecture; Laurent polynomial; Peano derivative; vector spaces and spans.}

%\begin{abstract}
%For each positive integer $n$, function $f$, and point $x$, the 1998 conjecture by Ghinchev, Guerragio, and Rocca states that the existence of the $n$-th Peano derivative $f_{(n)}(x)$ is equivalent to the existence of all $n(n+1)/2$ generalized Riemann derivatives,
%\[
%D_{k,-j}f(x)=\lim_{h\rightarrow 0}\frac 1{h^{k}}\sum_{i=0}^k(-1)^i\binom{k}{i}f(x+(k-i-j)h),
%\]
%for $j,k$ with $0\leq j<k\leq n$. A version of it for $n\geq 2$ replaces all $-j$ with $j$ and eliminates all $j=k-1$. Both the GGR conjecture and its version were recently proved by the authors using non-inductive proofs based on highly non-trivial combinatorial algorithms. This article provides a simple, inductive, algebraic proof of each of these theorems, based on a reduction to (Laurent) polynomials.
%\end{abstract}
\maketitle

\noindent Arguably, the most important application of higher order derivatives is Taylor's theorem, asserting that an $n$ times differentiable function $f$ at a point $c$ is approximated near $c$ by a polynomial $p$ with error $f(c+h)-p(c+h)=o(h^n)$ as $h\rightarrow 0$. It is also well known that Taylor's theorem provides only a sufficient condition for this approximation to happen, and all functions $f$ with this property are said to be $n$ times Peano differentiable at $c$.

\medskip
In 1998, Ginchev, Guerragio, and Rocca (GGR) conjectured the following result:
\begin{GGRconjecture}
When $n\geq 2$, the following two conditions,
\begin{enumerate}
\item[(i)] $f$ is $n-1$ times Peano differentiable at $c$ and
\item[(ii)] $\displaystyle \lim_{h\rightarrow 0}\frac 1{h^{n}}\sum_{j=0}^{n}(-1)^{n-j}\binom{n}{j}f(c+(j-k)h)$ exists for all $k$ with $0\leq k\leq n-1$,
\end{enumerate}
are sufficient to make $f$ an $n$ times differentiable function at $c$.
\end{GGRconjecture}

They proved the theorem by hand for $n=2,3,4$ in \cite{GGR}, and with the use of a computer they proved it for $n=5,6,7,8$ in \cite{GR}, leaving the rest as a conjecture. The GGR conjecture was recently proved in \cite{AC2} and is now a theorem. The result proved in \cite{AC2} is slightly stronger for $n$ odd, by eliminating the second condition for $k=0$. A variant of the conjecture, where the bounds for $k$ are replaced by $-(n-2)\leq k\leq 0$, is proved in \cite{ACF}. Article \cite{C} sheds some light  towards extending the GGR conjecture for $n=1$.

The original statement of the GGR conjecture is actually an equivalent version of the above, by the principle of mathematical induction: If $D_{n,k}$ denote the above limits, Ginchev, Guerragio, and Rocca conjectured that the $\binom {n+1}2$ limits $D_{m,k}$ for $0\leq k<m\leq n$ will be enough for $f$ to be $n$ times Peano differentiable at $c$. The limit $D_{n,0}=D_{n,0}f(c)$ is called the $n$-th Riemann derivative of $f$ at $c$, and $D_{n,n/2}f(c)$ is the $n$-th symmetric Riemann derivative of $f$ at $c$. Both of these derivatives were invented by Riemann in the mid 1800s, see \cite{R}. The Peano derivatives were invented by Peano in \cite{P} in 1892, and then developed greatly by de la Vall\'ee Poussin in \cite{dlVP}. For this reason, they are often referred to as de la Vall\'ee Poussin derivatives.

The purpose of this note is to provide a new, simple proof to both the GGR conjecture and its variant.
\[
\ast\quad\ast\quad\ast
\]

\medskip
Let $R_n(h)$ be the difference defined recursively by $R_1(h)=f(c+h)-f(c)$, and $R_n(h)=R_{n-1}(2h)-2^{n-1}R_{n-1}(h)$ for $n\geq 2$. Closed form formulas for generalizations of these differences are deduced in \cite{ACF1}; they involve the Gaussian or $q$-binomial coefficients, so they are $q$-analogues of the Riemann differences. Other $q$-analogues of Riemann differences are found in \cite{AC} and \cite{ACR}.

We will use the following 1936 result of Marcinkiewicz and Zygmund in \cite{MZ}.

\begin{theorem}\label{TT1}
Suppose $f$ is $n-1$ times Peano differentiable at $c$.
If $\,\lim_{h\rightarrow 0}R_n(h)/h^n$ exists, then $f$ is $n$ times Peano differentiable at $c$.
\end{theorem}

If we denote $\Delta_k(h)$ as the difference $\sum_{j=0}^{n}(-1)^{n-j}\binom{n}{j}f(c+(j-k)h)$ in the GGR limit condition, then this condition can be concisely written as $\lim_{h\rightarrow 0}\Delta_k(h)/h^n$ exists. It is also obvious that if all GGR limit conditions are met then the following linear combination
\[
\lim_{h\rightarrow 0}\frac {\sum_kc_k\Delta_k(s_kh)}{h^n}\text{ exists, where $c_k,s_k$ are arbitrary real constants.}
\]

Hence, if we can show that $R(h)$ is a linear combination $\sum_kc_k\Delta_k(s_kh)$, then the GGR conjecture follows from Theorem~\ref{TT1}.

That this is indeed the case will follow from the analogous result for polynomials, via the linear isomorphism,
$
\Delta(h):=\sum c_jf(c+b_jh)\mapsto d(t):=\sum c_jt^{b_j},
$
from the $\mathbb{R}$-space of all differences of $f$ at $c$ and $h$ with integer nodes (the $b_j$), to the $\mathbb{R}$-space $\mathbb{R}[t, t^{-1}]$ of all Laurent polynomials in indeterminate $t$ with real coefficients. In this way, (1) if $\Delta_k(h)=\sum_{j=0}^{n}(-1)^{n-j}\binom{n}{j}f(c+(j-k)h)$, then $d_k(t)=\sum_{j=0}^{n}(-1)^{n-j}\binom{n}{j}t^{j-k}=t^{-k}(t-1)^n$; (2) the polynomial corresponding to $\Delta_k(sh)$ is $d_k(t^s)$; and (3) 
if $r_n(t)$ is the polynomial that corresponds to $R_n(h)$ under this linear isomorphism, then its recursive definition is~$r_1(t)=t-1$, and $r_n(t)=r_{n-1}(t^2)-2^{n-1}r_{n-1}(t)$ for $n\geq 2$.

Based on these properties of the above linear isomorphism, we will be done by showing that the following result is true.

\begin{theorem}\label{TT2}
There are constants $c_k$ and $s_k$ such that $r_n(t)=\sum c_kd_k(t^{s_k})$.
\end{theorem}

Before proceeding with the proof of Theorem~\ref{TT2}, we need to make a clarification. Our solution to the theorem has the numbers $s_k$ non-negative integers instead of real numbers, so we can think of $s_k$ as $s$, with $s\geq 0$. Then the $c_k$ are viewed as $c_{k,s}$, for a more precisely indexed sum
$\sum c_kd_k(t^{s_k})=\sum_{s=0}^{\infty}\sum_kc_{k,s}d_k(t^s)=\sum_{s=0}^{\infty}\sum_kc_{k,s}t^{-sk}(t^s-1)^n$, where the ranges for $k$ will be clarified later, since these are different for different cases in the proof of the theorem. Summarizing, in order to prove Theorem~\ref{TT2}, it suffices to show that
\[r_n(t)\in V_n:=\mbox{span}\{t^{-sk}(t^s-1)^n\mid k=(0),1,\ldots,n-1,s=1,2,\ldots \},\]
where $(0)$ means that the value $0$ is taken only for $n$ even.
\[
\ast\quad\ast\quad\ast
\]

The proof of Theorem~\ref{TT2} is much shorter for the variant of the GGR conjecture than it is for the conjecture itself. For this reason, we deal with the variant first.

%%%%%%%%%%%%%%
\subsection*{Proof of Theorem~\ref{TT2} (Variant Case)} In this case, by replacing the index of summation~$k$ with $-k$, the range $-(n-2)\leq k\leq 0$ becomes $0\leq k\leq n-2$. In this way, $V_n$ becomes
\[
V_n=\mbox{span}\{t^{sk}(t^s-1)^n\mid k=0,\ldots,n-2;\;s=1,2,\ldots\}.
\]

The following lemma provides a new set of generators for the space $V_n$.

\begin{lemma}\label{L8}
$V_n={\rm span}\{(t^s-1)^{n+k}\mid k=0,1,\ldots ,n-2,\; s=1,2,\ldots  \}$.
\end{lemma}

\begin{proof}
It suffices to show the following equality of subspaces:
\[
\mbox{span}\{t^{k}(t-1)^n\mid k=0,1,\ldots ,n-2\}={\rm span}\{(t-1)^{n+k}\mid k=0,1,\ldots ,n-2 \}.
\]
Indeed, this is the result of multiplying by $(t-1)^n$ both sides of the obvious equation
$\mbox{span}\{1,t,t^2,\ldots,t^{n-2}\}={\rm span}\{1,t-1,(t-1)^2,\ldots ,(t-1)^{n-2} \}$.
\end{proof}

We are now ready to proceed with the proof of Theorem~\ref{TT2} in its variant case.

\begin{proof}[Proof of Theorem~\ref{TT2} (Variant Case)]
Induct on $n$. When $n=2$, $r_2=(t-1)^2$ is clearly in~$V_2$. Suppose $r_n\in V_n$, for some $n$, $n\geq 2$, and prove the same property for $n+1$. By Lemma~\ref{L8},~$r_n$ is a linear combination of polynomials of the form $(t^s-1)^{n+k}$, where~$s$ is a positive integer and $k=0,1,\ldots ,n-2$.
By the recursion, $r_{n+1}(t)=r_n(t^2)-2^nr_n(t)$ will be a linear combination of polynomials
\[
(t^{2s}-1)^{n+k}-2^n(t^s-1)^{n+k}\text{, for various $k$ and $s$}.
\]
By Lemma~\ref{L8}, these polynomials belong to $V_{n+1}$ in all cases, except for $k=0$,~when
\[
(t^{2s}-1)^{n}-2^n(t^s-1)^{n}=(t^s-1)^n((t^s+1)^n-2^n)=(t^s-1)^{n+1}p(t^s),
\]
where $p$ is a polynomial in $t$ of degree $n-1$, so that $(t^s-1)^{n+1}p(t^s)$ belongs to the 
subspace $\mbox{span}\{ (t^{s}-1)^{n+1},(t^{s}-1)^{n+2},\ldots ,(t^{s}-1)^{2n}\}$ of $V_{n+1}$. 
\end{proof}

%%%%%%%
\subsection*{Proof of Theorem~\ref{TT2} (GGR Case)} The GGR case in the proof of Theorem~\ref{TT2} is similar to the variant case.
The proof of the inductive step is now split further into two subcases, $n$ odd and $n$ even.
In both cases, following the refined result of the GGR Theorem from \cite{AC2},
\[
V_n=\mbox{span}\{t^{-sk}(t^s-1)^n\mid k=(0),1,\ldots,n-1;\;s=1,2,\ldots\},
\]
where $(0)$ means that $0$ is taken for $n$ even, and not taken for $n$ odd. 
More explicitly, this~is
\[
V_n=\begin{cases}
\mbox{span}\{ (t-1)^n,t^{-1}(t-1)^n,\ldots ,t^{-(n-1)}(t-1)^n,\ldots \},&n\text{ even},\\
\mbox{span}\{ t^{-1}(t-1)^n,t^{-2}(t-1)^n,\ldots ,t^{-(n-1)}(t-1)^n,\ldots \},&n\text{ odd},
\end{cases}
\]
where the last dots in both cases mean that the generating set also includes the previously listed generators evaluated at $t^s$, for all $s$ at least 2.
Let $W_n$ be the subspace of $V_n$ spanned by all generators with $s=1$. Then $W_n$ has the expression
\[
W_n=\begin{cases}
\mbox{span}\{ (t-1)^n,t^{-1}(t-1)^n,\ldots ,t^{-(n-1)}(t-1)^n\},&n\text{ even},\\
\mbox{span}\{ t^{-1}(t-1)^n,t^{-2}(t-1)^n,\ldots ,t^{-(n-1)}(t-1)^n\},&n\text{ odd}.
\end{cases}
\]
The following lemma provides new sets of generators for $W_n$ in both parity cases.
\begin{lemma}\label{L9}
With the above notation,
\[
W_n=\begin{cases}
{\rm span}\{ t^{-n/2}(t-1)^n,t^{-1}(t-1)^{n+1},\ldots ,t^{-(n-1)}(t-1)^{n+1}\},&n\text{ \rm even},\\
{\rm span}\{ t^{-(n-1)/2}(t-1)^n,t^{-2}(t-1)^{n+1},\ldots ,t^{-(n-1)}(t-1)^{n+1}\},&n\text{ \rm odd}.
\end{cases}
\]
\end{lemma}

\begin{proof}
When $n$ is even, the result follows from $t^{-n/2}(t-1)^n$ being one of the generators in the definition of $W_n$, and $t^{-k}(t-1)^{n+1}=t^{-(k-1)}(t-1)^{n}-t^{-k}(t-1)^{n}$, for each $k=1,\ldots ,n-1$. The case when $n$ is odd has a similar proof.
\end{proof}

We are now ready to prove Theorem~\ref{TT2} in the GGR case.

\begin{proof}[Proof of Theorem~\ref{TT2} (Case GGR)]
Induct on $n$. When $n=1$, $r_1=t-1\in V_1$.
We assume that $r_n\in V_n$
and prove that $r_{n+1}\in V_{n+1}$ in two possible cases:

\textbf{Case 1.} When $n$ is even,
\[
\qquad V_{n+1}=\mbox{span}\{ t^{-1}(t-1)^{n+1},t^{-2}(t-1)^{n+1},\ldots ,t^{-n}(t-1)^{n+1},\ldots \}.
\]
The inductive hypothesis and Lemma~\ref{L9} imply that $r_n(x)$ is a linear combination of
\[
t^{-n/2}(t-1)^n\text{ and }t^{-k}(t-1)^{n+1}\text{, for }k=1,\ldots,n-1,
\]
and their evaluations at $t^s$, for $s$ at least~2. Then $r_{n+1}(t)=r_n(t^2)-2^nr_n(t)$ is a linear combination of two kinds of polynomials and their evaluations at $t^s$, for~$s\geq 2$. The first kind of polynomial has the form $t^{-n}(t^2-1)^n-2^nt^{-n/2}(t-1)^n$
\[
=t^{-n}(t-1)^n((t+1)^n-2^nt^{n/2})=t^{-n}(t-1)^{n+1}p(t),
\]
where $p(t)$ is a polynomial degree $n-1$, hence the whole expression lives inside of
\[
(t-1)^{n+1}{\rm span}\{t^{-1},t^{-2},\ldots ,t^{-n}\},
\]
a subspace of $V_{n+1}$. The polynomials of the second kind are polynomials of the~form
\[
t^{-2k}(t^2-1)^{n+1}-2^nt^{-k}(t-1)^{n+1}\text{, for }k=1,\ldots,n-1.
\]
Their second term is a scalar multiple of a generator of $V_{n+1}$, while their first term is an $(s=2)$-dilation of the same generator, so all polynomials of the second kind also belong to $V_{n+1}$. We conclude that $r_{n+1}\in V_{n+1}$, as needed.

\textbf{Case 2.} When $n$ is odd,
\[
\qquad V_{n+1}=\mbox{span}\{ (t-1)^{n+1},t^{-1}(t-1)^{n+1},\ldots ,t^{-n}(t-1)^{n+1},\ldots \}.
\]
The inductive hypothesis and Lemma~\ref{L9} imply that $r_n(x)$ is a linear combination of
\[
t^{-(n-1)/2}(t-1)^n\text{ and }t^{-k}(t-1)^{n+1}\text{, for }k=2,\ldots,n-1,
\]
and their evaluations at $t^s$, for $s$ at least~2. Then $r_{n+1}(t)=r_n(t^2)-2^nr_n(t)$ is a linear combination of two kinds of polynomials and their evaluations at $t^s$, for~$s\geq 2$. The first kind of polynomial is of the form $t^{-(n-1)}(t^2-1)^n-2^nt^{-(n-1)/2}(t-1)^n$
\[
=t^{-(n-1)}(t-1)^n((t+1)^n-2^nt^{(n-1)/2})=t^{-(n-1)}(t-1)^{n+1}p(t),
\]
where $p(t)$ is a polynomial degree $n-1$, hence the above expression lives inside~of
\[
(t-1)^{n+1}{\rm span}\{1,t^{-1},t^{-2},\ldots ,t^{-(n-1)}\},
\]
a subspace of $V_{n+1}$. The polynomials of the second kind are polynomials of the~form
\[
t^{-2k}(t^2-1)^{n+1}-2^nt^{-k}(t-1)^{n+1}\text{, for }k=2,\ldots,n-1.
\]
Their second term is a scalar multiple of a generator of $V_{n+1}$, while their first term is an $(s=2)$-dilation of the same generator, so all polynomials of the second kind also belong to $V_{n+1}$. We conclude that $r_{n+1}\in V_{n+1}$ in this case as well.
\end{proof}
%%%%%%%%%%%%%%

%%%%%%%%%%%%%%%%%%%%%%%%%%%%%%%%%%%%%%%%%%%%%%%%%%%%%%}
\bibliographystyle{amsplain}

\end{document}